\documentclass{article}
\usepackage{amsmath,amsthm,amsfonts,amscd,amssymb,eucal,latexsym}

\begin{document}

\title{Ergodic actions of $S_\mu U(2)$ on $C^*$--algebras from $II_1$ 
subfactors}
\author{Claudia Pinzari$^1$, John E. Roberts$^2$\\ \\
$^1\,$Dipartimento di Matematica, Universit\`a di Roma ``La Sapienza''\\
00185--Roma, Italy\\
$^2\,$Dipartimento di Matematica, Universit\`a di Roma ``Tor Vergata'',\\
00133 Roma, Italy}
\date{}
\maketitle

\begin{abstract}
To  a proper inclusion $N\subset M$ of $II_1$ 
factors of finite Jones  index $[M:N]$,
we associate an ergodic $C^*$--action of the quantum group 
$S_{\mu}U(2)$  
(or more generally of certain groups $A_o(F)$). The
higher relative commutant
$N'\cap M_{r-1}$ can be 
identified with the  spectral space of the $r$-th tensor 
power $u^{\otimes r}$ of the defining representation of the quantum group.
The index and the deformation parameter are related by $-1\leq\mu<0$ and 
$[M:N]=|\mu+\mu^{-1}|$.

This ergodic action may be thought of as a virtual subgroup  of
 $S_\mu U(2)$ in the sense of Mackey arising from  the tensor 
 category generated by the   
$N$--bimodule 
${}_N M_N$.
 $\mu$ is negative as
${}_NM_N$ is a real bimodule.
\medskip



\noindent{\it Subj. Class.}: quantum groups, noncommutative topology and geometry, dynamical systems

\end{abstract}

\begin{section} {Introduction}

Ergodic theory of  compact quantum groups 
on unital 
$C^*$--algebras is the noncommutative analogue of 
the theory of 
homogeneous spaces over compact groups. 
Several examples are known, showing that this theory 
exhibits 
new aspects with respect to the classical case.

The first interesting examples, 
the quantum spheres, have been constructed by Podle\'s
\cite{Podles1}, and   the theory has been continued by Boca  \cite{Boca}.
While in the commutative case, homogeneous 
spaces are quotients by closed subgroups, some of Podle\'s examples show that 
is no longer true in the 
noncommutative situation. Hence, an ergodic $C^*$--action of a compact quantum 
group may be regarded as a 
virtual subgroup, introduced by
Mackey   in classical ergodic theory  
\cite{Mackey1}, \cite{Mackey2},  but in a compact and 
noncommutative setting. 

Whereas the invariant state of a $C^*$--algebra carrying an ergodic action 
of a compact group is necessarily tracial \cite{HLS}, this is not the case of the
ergodic actions of the van Daele-Wang quantum groups 
on type   $III_\lambda$ factors constructed 
by \cite{Wang2}. 
Further examples of ergodic actions have been given in \cite{BRV} starting from  
 tensor equivalences of representation categories of compact  quantum 
groups.

 An axiomatization of the spectral functor 
associated to an ergodic $C^*$--action has been given
in \cite{PR}, where
such functors 
were called {\it quasitensor}.

The aim of this paper is to
establish a connection between the ergodic theory of compact quantum groups
on unital $C^*$--algebras and the Jones theory of subfactors.

In his pioneering work, Jones showed that the theory of 
$II_1$ subfactors  with finite index can be understood as a generalized 
group theory, exhibiting  beautiful
 representations of the infinite braid group 
 factoring through the Temperley-Lieb relations.
A well known consequence is that
Ocneanu's bimodule category  associated with  
 the inclusion, contains, roughly speaking, a deformation of  the
 representation category 
 of $SU(2)$ generated by the flip map \cite{Jones1}, \cite{Jones2}, \cite{GHJ}. 
But as this deformed category is not associated with representations of the symmetric groups, 
this bimodule category is not described by a group, nor even by
a quantum group, in general.

Our main point is that  
a finite index inclusion of $II_1$ factors always 
gives rise to an ergodic action of 
 $S_\mu U(2)$ for suitable  negative values of $\mu$.
To understand why the values are negative note that 
the defining representation of 
$SU(2)$ is selfconjugate and  {\it pseudoreal} 
whilst $M$, as a $N$--bimodule, is selfconjugate and {\it real}, 
However, the defining representation of the quantum group $S_\mu 
U(2)$  of Woronowicz \cite{WoronowiczTK} is real provided the 
deformation parameter is {\it negative} (see, e.g., \cite{P1}).

More precisely, we canonically associate  to any proper inclusion 
of $II_1$ factors with finite index, an ergodic action of the van 
Daele--Wang quantum groups $A_o(F)$ on unital $C^*$--algebras,
 where the spectral spaces of the action are the higher relative 
commutants $N'\cap M_{r}$. 
 $F$ is subject to the conditions $F\overline{F}=I$ and 
$\text{Tr}(F^*F)=[M:N]$.
By a result of Banica \cite{Banica0} (see also \cite{PR}),
whilst  the quantum group is not unique, its 
representation category is determined by the above conditions up to a 
tensor isomorphism. 
In particular, choosing $F$ of minimal rank,
yields  ergodic actions of $S_\mu U(2)$ where the index and the 
deformation parameter are related by 
$-1\leq\mu<0$ and $|\mu+\mu^{-1}|=[M:N]$. On the other hand, when the index 
is an integer, the identity matrix of rank $[M:N]$ is a natural choice, providing an 
ergodic action of the Kac type quantum group $A_o([M:N])$. In this case, 
we show that the invariant state of the associated $C^*$--algebra is tracial.

Taking Wenzl's work \cite{Wenzl1}, \cite{Wenzl2}, \cite{Wenzl3} on 
constructing  subfactors from algebraic quantum groups  into account   
and passing through subfactors we connect quantum groups at roots of 
unity to ergodic actions of compact quantum groups whose quantum 
dimension depends on roots of unity.
Thus quantum groups at roots of unity would seem to be  
virtual quantum subgroups of compact quantum groups.

However,  the examples of Asaeda--Haagerup \cite{AH}
show that not all subfactors are associated with quantum groups. 
Correspondingly, we give examples
of ergodic actions of $S_\mu U(2)$, for $\mu<0$, not arising from 
quantum group constructions. This is a novelty compared to the group case,
as Wassermann has shown that all ergodic actions of $SU(2)$ arise from 
closed subgroups and their irreducible projective representations via 
induction constructions \cite{Wassermann3}.

Our approach uses the duality theorem of 
\cite{PR}. We start from the remark  that the most important 
axiom of that theorem (see $(2.4)$)
is in analogy with the  {\it commuting
square condition} of \cite{PopaMaximal}, \cite{Pi-Po}, which plays a key role 
in Jones's theory of subfactors. This analogy
raises the question of 
whether this theory yields quasitensor functors and hence ergodic 
$C^*$--actions of compact quantum groups.

The first step of our construction is to show that 
if $\mu$ is suitably chosen,  the full tensor 
$C^*$--subcategory 
of Rep(S$_\mu$U(2)) generated by the fundamental representation embeds
into the $SU(2)$--like category contained in the bimodule tensor category 
 generated by  ${}_NM_N$.

In particular,  from the pivotal result of Jones 
restricting the values of the index
\cite{Jones1},  only the  values $4\cos^2\pi/m$, $m\geq 4$  and $\geq 4$ of the 
quantum dimension of $S_\mu U(2)$ can possibly arise in our examples.

 Because of finiteness of the factors in question, the quantum multiplicity of  
 $u^{\otimes r}$ in the ergodic action, where $u$ is the defining 
representation of  $A_o(F)$,   takes its lower bound, the integral multiplicity,  
in turn given by the dimension of the higher relative 
commutant $N'\cap M_{r-1}$. 
Moreover the spectral space corresponding to $u^{\otimes r}$ with its Hilbert space structure can be identified with the 
higher relative commutant $N'\cap M_{r-1}$, with inner product coming from the Markov trace.

Furthermore these ergodic actions are not, in general,
embedable into the translation action. 
We shall show  this  for  non-integral 
values of the index of an extremal and amenable inclusion. 
The proof relies on Popa's work \cite{Popa}. 

We show in \cite{PR2}  that the results of this paper 
have analogues when subfactors are replaced by tensor 
$C^*$--categories with conjugation. 
 The groups involved are $A_u(F)$ and $A_o(F)$.
Such categories, arise, in particular, in the algebraic 
approach to QFT where they are in addition endowed 
with a unitary braided symmetry,
(see, e.g. \cite{Haag}).

In a subsequent paper, we will adopt Mackey's point of view
that a non-transitive ergodic action can be viewed 
as a virtual subgroup which should thus exhibit 
typical properties  of a closed subgroup.
We will develop a theory of induction and restriction for 
representations of these virtual subgroups
 \cite{DR2}.
\medskip

The paper is organized as follows. In section 2 we recall the main 
invariants of ergodic $C^*$--actions and the 
duality theorem of \cite{PR}.

In section 3 we recall Ocneanu's bimodules associated with an 
inclusion of $II_1$ factors with finite Jones index and we show  that 
extremal and amenable inclusions in the sense of \cite{Pi-Po}, \cite{Popa}
give rise to categories which, for non-integral values of the index, are 
not 
embedable into the category of Hilbert spaces. When they are, only 
quantum groups with coinvolutive coinverses appear.
Hence Ocneanu's categories associated to amenable inclusions are, in this 
respect, rather different 
from the representation categories of compact quantum groups, which are 
embedded in Hilbert spaces by construction, but often give rise to 
non-amenable inclusions.

In section 4 we state our main results: the 
inclusion gives rise to ergodic
$C^*$--actions of the van Daele-Wang compact quantum group $A_o(F)$ 
associated with an invertible matrix
$F\in M_n({\mathbb C})$.

Section 5 is devoted to the proofs of the results.
We conclude the paper  with the necessary computations to yield a 
presentation  of the dense
 $^*$--subalgebra of 
spectral elements  by generators and 
relations, in terms of the higher 
relative commutants $N'\cap M_r$.

 \end{section}

\begin{section}{Preliminaries on ergodic $C^*$--actions}

Consider a unital $C^*$--algebra ${\cal C}$ and a compact quantum group 
$G=({\cal Q},\Delta)$, with coproduct
$\Delta:{\cal Q}\to{\cal Q}\otimes 
{\cal Q}$ 
\cite{WLesHouches}.
An {\it action} of $G$ on ${\cal C}$ is a unital $^*$--homomorphism
$\delta:{\cal C}\to{\cal C}\otimes{\cal Q}$ satisfying the group 
representation property:
$$\iota\otimes\Delta\circ\delta=\delta\otimes\iota\circ\delta,$$
and the nondegeneracy property requiring  $\delta({\cal 
C})I\otimes{\cal 
Q}$ to be dense in ${\cal C}\otimes{\cal Q}$. The {\it spectrum} of 
$\delta$,
$\text{sp}({\delta})$, is defined  to be the set of all 
representations $u$ of $G$ for which there is a faithful linear map
$T: H_u\to{\cal C}$ intertwining the representation $u$ with the action 
$\delta$:
$$\delta\circ T=T\otimes\iota\circ u.$$
In other words, if $u_{ij}$ are the 
coefficients of $u$ corresponding to some orthonormal basis of $H_u$,
we are requiring the existence of linearly independent elements
$c_1,\dots, c_d\in{\cal C}$, with $d$ the dimension of $u$, transforming 
like $u$ under the action:
$\delta(c_i):=\sum_j c_j\otimes u_{ji}$. The linear span of all the 
$c_i$'s, 
denoted ${\cal C}_{\text{sp}}$,
as $u$ varies in the spectrum, is a dense $^*$--subalgebra of ${\cal C}$
\cite{Podles}.
 
The action $\delta$ is called {\it ergodic} if the fixed point algebra
$${\cal C}^\delta=\{c\in{\cal C}:\delta(c)=c\otimes I\}$$ reduces to the 
complex numbers: ${\cal C}^\delta={\mathbb C}I.$
The simplest example of an ergodic action is the translation action of $G$ 
on
${\cal C}={\cal Q}$ with $\delta=\Delta$.

If an action $\delta $ is ergodic,   
spectral 
multiplets can be organized to form 
Hilbert spaces.
In fact, for any representation $u$, consider the space
$$L_u:=\{T:H_u\to{\cal C}, \delta\circ T =T\otimes\iota\circ u\}.$$
If $S,T\in L_u$, 
$<S,T>:=\sum_iT(\psi_i)S(\psi_i)^*$, with $(\psi_i)$ an orthonormal basis 
of $H_u$, is an element of the 
fixed point 
algebra
${\cal C}^\delta$, and therefore a complex number.
$L_u$ is known to be finite dimensional, and therefore  a 
Hilbert space with the above inner product.
This Hilbert space is nonzero precisely when $u$ contains a 
subrepresentation $v\in\text{sp}(\delta)$. In particular, for an 
irreducible $u$, the conditions $u\in\text{sp}(\delta)$ and $L_u\neq 0$ 
are equivalent. The dimension of $L_u$ is called the {\it multiplicity}
of $u$ and denoted $\text{mult}(u)$.

The complex conjugate vector space $\overline{L_u}$,
 endowed with the conjugate inner product 
$$<\overline{S},\overline{T}>:=<T,S>=\sum_i S(\psi_i)T(\psi_i)^*,$$
is called the {\it spectral space} associated with $u$.

If for example  $\delta$ is the translation action on ${\cal Q}$, any
$\psi\in H_u$ defines an element of $L_u$ by
$$T_\psi(\psi'):=\psi^*\otimes Iu(\psi').$$
Hence the spectral space $\overline{L_u}$ 
can be identified with $H_u$ through the unitary map
$$\psi\in H_u\to \overline{T_\psi}\in \overline{L_u}.$$

If $j:H_u\to H_{\overline{u}}$ is an antilinear invertible defining a conjugate (unitary) representation
 $\overline{u}$ of
 $u$ by
$\overline{u}_{\phi,j\psi}=(u_{j^*\phi,\psi})^*$, then
there is an associated 
antilinear 
$J:L_u\to L_{\overline{u}}$ by $J(T)(\phi):=T(j^{-1}(\phi))^*$
with inverse $J^{-1}:L_{\overline{u}}\to L_u$ given
by $J^{-1}(S)(\psi)=S(j(\psi))^*$. If $u$ is irreducible,
 the {\it quantum multiplicity} $m(u)$ of $u$ is defined by
$m(u)^2:=\text{Trace}(JJ^*)\text{Trace}((JJ^*)^{-1})$ \cite{BRV}.
One has:
$$\text{mult}(u)\leq m(u)\leq d(u),$$
an inequality which strenghtens the inequality $\text{mult}(u)\leq d(u)$,
with $d(u)$ the quantum dimension of $u$,
previously obtained by Boca \cite{Boca}, in turn generalizing HLS theorem 
\cite{HLS}
$\text{mult}(u)\leq \text{dim}(u)$ in the group case.
If $u$ is reducible, we define $m(u)$ as  the infimum of all 
the above trace values, when $j$ ranges over all possible solutions of the 
conjugate equations. Then the  inequality 
$$\text{dim}(L_u)\leq m(u)\leq 
d(u)\eqno(2.1)$$
now  holds for all representations $u$.
Notice that $m(u)$ takes the smallest possible value $\text{dim}(L_u)$ 
precisely when for some $j$ the associated $J$ is a scalar multiple of an
antiunitary. The actions we shall construct in this paper satisfy this 
property. 
Examples of ergodic actions of $S_\mu U(2)$ where 
$\text{dim}(u)<\text{mult}(u)<m(u)=d(u)$ have been constructed in 
\cite{BRV}.

The study of  categorical aspects of ergodic $C^*$--actions of compact 
quantum groups has been 
developed in 
 generality  in \cite{PR}, 
where a spectral characterization has been obtained.
It has been shown that the {\it spectral functor} of an ergodic 
$G$--action
is a dual object, in the sense that 
the 
$G$--action on the maximal completion of 
 ${\cal C}_{\text{sp}}$ can be reconstructed from it. Furthermore, 
spectral functors of ergodic $C^*$--actions of compact quantum groups 
are  
characterized among all $^*$--functors from $\text{Rep}(G)$ to the 
category of Hilbert spaces, by the property of being {\it quasitensor}.
More in detail, the spectral functor $$\overline{L}:\text{Rep}(G)\to{\cal 
H}$$ 
associated with an ergodic $G$--action is the 
functor from the category  of representations of $G$
to the category ${\cal H}$ of Hilbert spaces. 
This functor is defined as $u\to\overline{L_u}$ 
on objects and as follows on arrows.

If $A\in(u,v)$ and $T\in L_v$ then $T\circ A: H_u\to{\cal C}$
lies in $L_u$.
Hence 
if we identify 
$\overline{L_u}$ canonically with the dual vector space of $L_u$, any 
arrow
$A\in(u,v)$ in $\text{Rep}(G)$ induces a linear map
$\overline{L}_A\in(\overline{L_u},\overline{L_v})$
by 
$$\overline{L}_A:\varphi\in\overline{L}_u\to (T\in L_v\to \varphi(T\circ 
A))\in\overline{L_v}.$$
Taking into account the tensor $C^*$--category structure of 
$\text{Rep}(G)$ and ${\cal H}$ one can see that $\overline{L}$ becomes a 
$^*$--functor, but {\it not} a tensor $^*$--functor, in general.

As far as the tensor structure of $\overline{L}$ is concerned, for 
$u,v\in\text{Rep}(G)$, the 
tensor product Hilbert 
space 
$\overline{L_u}\otimes\overline{L_v}$ 
is in general just a subspace
of $\overline{L_{u\otimes v}}$, in the sense that there is a 
natural isometric inclusion
$$\lambda_{u,v}:\overline{L_u}\otimes\overline{L_v}\to\overline{L_{u\otimes 
v}}$$
identifying any simple tensor $\overline{S}\otimes\overline{T}$ with
the complex conjugate of the element of $L_{u\otimes v}$ defined by
$$\psi\otimes\phi\in H_u\otimes H_v\to S(\psi)T(\phi).$$
The dual of the action is the pair $(\overline{L}, \lambda)$ 
consisting of the 
functor $\overline{L}$ and  all the inclusions $\lambda_{u,v}$.

The maximal $C^*$--completion of ${\cal 
C}_{\text{sp}}$, together with the extended $G$--action, can be 
reconstructed 
from the dual $(\overline{L}, \lambda)$.

The main result of \cite{PR} is an axiomatization of the 
set of all duals of ergodic 
actions 
$(\overline{L},\lambda)$ among all
$^*$--functors
$${\tau}:\text{Rep}(G)\to{\cal H}$$
endowed  with isometries $\tilde{\tau}_{u,v}:{\tau}_u\otimes\tau_v\to{\tau}_{u\otimes v}$. 
All pairs $(\tau, \tilde\tau)$ satisfying  
 properties $(3.1)$--$(3.6)$  in \cite{PR} have been shown to arise as 
the dual of an ergodic $G$--action. Such functors were called 
quasitensor. In \cite{PR2} the following equivalent simpler
axiomatization has been derived:
$${\tau}_\iota=\iota,\eqno(2.2)$$
$$\tilde\tau_{u,\iota}=\tilde\tau_{\iota, u}=1_{\tau_u},\eqno(2.3)$$
$$\tilde\tau_{u,v\otimes w}^*\circ\tilde\tau_{u\otimes v,w}=
1_{\tau_u}\otimes\tilde\tau_{v,w}
\circ{\tilde\tau_{u,v}}^*\otimes 1_{\tau_w}\eqno(2.4)$$
$$\tau({S\otimes T})\circ
\tilde\tau_{u,v}=\tilde\tau_{u',v'}\circ\tau(S)\otimes\tau(T),\eqno(2.5)$$
 for any other pair of objects $u'$, $v'$ and arrows
$S\in(u, u')$,
$T\in(v,v')$.
In particular, a
{\it tensor} functor $\tau$ is quasitensor with 
$\tilde\tau_{u,v}:=1_{\tau_u\otimes\tau_v}$.
An ergodic $C^*$--action of $G$ on a unital $C^*$--algebra ${\cal C}$ can 
be constructed by duality from a 
quasitensor $^*$--functor $(\tau, \tilde\tau):\text{Rep}(G)\to{\cal H}$. Once 
the ergodic action has been  constructed, the pair $(\tau, \tilde\tau)$ can be 
identified
with the 
dual object of 
that action.
\medskip

\end{section}

\begin{section} {Ocneanu's 
 category from  a $II_1$ inclusion}

Consider an inclusion of $II_1$ factors $N\subset M$ of finite Jones index
$[M:N]$ \cite{Jones1} and denote 
the trace-preserving conditional expectation by
$E:M\to N$.
Let 
$$N\subset M\subset M_1\subset M_2\dots$$
be the Jones tower of $II_1$ factors. We denote 
the $r$-th Jones projection derived from the trace-preserving conditional
expectation $E_{r-1}: M_{r-1}\to M_{r-2}$
by $e_r\in L^2(M_{r-1})$, where we have set
$M_0=M$, $M_{-1}=N$, $E_0=E$. 
Recall that 
$M_r:=M_{r-1}e_r M_{r-1}$ and that $E_r(e_r)=[M:N]^{-1}$.
Also recall that the algebras $N'\cap M_r$, usually called the higher 
relative commutants, are
finite dimensional and 
$$\text{dim}(N'\cap M_{r-1})\leq [M:N]^r,\quad r\geq 0.$$
The main relations are the following. For $r\geq-1$,
$$[M_r, e_{r+2}]=0,\eqno(3.1)$$
$$e_{r+1}me_{r+1}=E_{r}(m)e_{r+1},\quad m\in M_r,\eqno(3.2)$$
implying the Jones projection relations:
$$e_ie_j=e_je_i,\quad |i-j|\geq2,\eqno(3.3)$$
$$e_ie_{j\pm1}e_i=[M:N]^{-1} e_i.\eqno(3.4)$$ 
We  review the well known construction of the Jones tower, $M_r$, 
$r\geq0$,  in terms 
of  Ocneanu's bimodules  \cite{Ocneanu88}, \cite{Pi-Po}, 
\cite{Pimsner-Popa}, \cite{Bisch}.
Regard $M$ as a right Hilbert 
$N$--module  with $N$--valued inner product
$$<m,m'>:=E(m^*m'),\quad m,m'\in M.$$
Since the index is finite, $M$ is finitely generated over $N$. Left 
multiplication on $M$ by elements of $N$ makes $M$ into a Hilbert bimodule
in the sense considered in  \cite{Pimsner}. Therefore we can take tensor 
powers of $M$
over $N$ and get further Hilbert bimodules. When no confusion arises, this tensor product will be simply denoted by $\otimes$.
As $N$--bimodules:
$$M\otimes M\simeq M_1.$$
Therefore
 iteratively, for $r=1,2,\dots$, 
$$M_{r-1}\simeq M^{\otimes r}.\eqno(3.5)$$
Consider the category ${\cal T}_M$ with objects 
the $N$--bimodules $M^{\otimes r}$, $r\geq0$, and arrows the 
bimodule 
mappings. This is 
a tensor $C^*$--category in the sense of \cite{DRInventiones}, with tensor 
product structure on arrows 
naturally induced
by the tensor product of Hilbert bimodules.
The arrow space 
$(\iota,M^{\otimes r})$ can be identified, as a vector space, with $N'\cap 
M_{r-1}.$
In particular,
$(\iota, M)\simeq N'\cap M\neq0$ and  is one-dimensional precisely 
when the inclusion is irreducible. These observations allow 
to show that
${\cal T}_{M}$, as a tensor $C^*$--category, is determined by an 
isometry $S\in(\iota, M)$ and 
its subcategory of arrow spaces 
$(M^{\otimes r}, M^{\otimes r})$ 
with same source and range objects.

Identifying the $N$--bimodules
$M^{\otimes r}\simeq M_{r-1}$,
the algebra 
${\cal L}_N(M^{\otimes r})$
 of right $N$--module maps
can be identified with Jones's basic construction associated with the 
inclusion
$N\subset M_{r-1}$, which, in turn, can be identified 
with
$M_{2r-1}$.
 The Jones projection of this inclusion
is given by $$f_{r-1}:=[M:N]^{r(r-1)/2}(e_r\dots e_1)(e_{r+1}\dots 
e_2)\dots(e_{2r-1}\dots e_r)$$ \cite{Pimsner-Popa}. Therefore as $C^*$--algebras,
$$(M^{\otimes r}, M^{\otimes r})\simeq N'\cap M_{2r-1}.$$
Only  terms of the Jones tower with odd indices appear
as we started with the bimodule ${}_N M_N$ rather than $\sigma:={}_M 
M_N$ or $\overline{\sigma}:={}_N M_M$. One has 
${}_N M_N\simeq\overline{\sigma}\otimes_M \sigma$.
It is well known and easy to check that tensoring on the {\it right} by 
$1_M$, 
namely
$T\in(M^{\otimes r}, M^{\otimes r})\to T\otimes 1_M\in(M^{\otimes r+1}, 
M^{\otimes r+1})$ 
corresponds to the natural inclusion $N'\cap 
M_{2r-1}\subset N'\cap M_{2r+1}$, whereas tensoring on the {\it left}
by $1_M$
corresponds to Ocneanu's {\it canonical shift} 
$\Gamma:N'\cap M_{2r-1}\to N'\cap M_{2r+1}$ \cite{Ocneanu88}.

We next give a   result showing that 
 inclusions of factors often
provide examples of categories ${\cal T}_M$ that cannot arise from 
representations of compact quantum 
groups.

In nongeneric cases, where ${\cal T}_M$ is embedable, it must generate 
the category of representations of a compact quantum group  with involutive 
coinverse. 
\medskip

\noindent{\bf 3.1 Proposition} {\sl Let $N\subset M$ be a finite
index, extremal and  amenable
inclusion of $II_1$ factors.
If $[M:N]$ is not an integer then
the tensor $C^*$--category ${\cal T}_M$ cannot be embedded into the
category of Hilbert spaces. Conversely, if $[M:N]$ is an integer
and if ${\cal T}_M$ is embedable then the Hilbert space
corresponding to ${}_NM_N$ has dimension $[M:N]$ and ${\cal T}_M$
is isomorphic to the representation category of a compact quantum subgroup
of the compact Kac quantum group  $A_o(I_{[M:N]})$ where $I_{[M:N]}$ 
is the
identity
matrix of size $[M:N]$.}\medskip

\noindent{\it Proof}
Popa  shows that, under the amenability assumption,
$$[M:N]=\lim_r\text{dim}(N'\cap M_r)^{1/r}$$
(see Theorem 4.4.1(3) in \cite{Popa}).
Therefore  if
 $[M:N]$ is not an integer, for sufficiently large $r$,
$$\text{dim}(\iota,M^{\otimes r})=\text{dim}(N'\cap M_{r-1})>i^r$$
where $i$ is the integral part of the index.
If ${\cal T}_M$ were embedable and if $n$ denotes
the dimension of the Hilbert space corresponding to $M$ then   
we must have
$n\leq [M:N]$. Hence $n\leq i$.
Furthermore for all $r$ we should also have
$$\text{dim}(\iota, M^{\otimes r})\leq n^r\leq i^r,$$
a contradiction.
Conversely, assume that $[M:N]=i$ and ${\cal T}_M$ embedable. If the
Hilbert space corresponding to $M$ had dimension  $n<i$ then
$$i=\lim_r\text{dim}(\iota,M^{\otimes r})^{1/r}\leq n<i,$$
again a contradiction. Therefore the dimension of the Hilbert space $H$
corresponding to $M$ is uniquely determined by the index. The
intertwiner $R$ needs then to correspond to an element $R'\in H\otimes H$
 with $\|R'\|^2=[M:N]=\text{dim}(H)$. We will show later that $M$
is a real object in the sense recalled at the beginning of Sect.\ 5 (see 
Theorem 5.2), hence $R'$ 
makes $H$ into  a 
real
object with
$\|R'\|^2=\text{dim}(H)$. This equality is
possible only if $R'=\sum_i e_i\otimes e_i$ for some orthonormal basis of
$H$, and the proof is complete. 
\medskip

\noindent{\it Remark} A notion of amenability for an object $\rho$ of a 
tensor $C^*$--category is introduced in \cite{LR}. It implies 
$d(\rho)=\lim_r\text{dim}(\rho^r,\rho^r)^{1/2r}$. As here, it is shown 
that ${\cal T}_\rho$ is not embedable unless $d(\rho)$ is an integer 
and that when it is then $d(\rho)=\text{dim} H(\rho)$.
\medskip

We illustrate the previous proposition with some known examples.
\medskip

\noindent a) {\it Fixed point and crossed products inclusions.}
The basic examples of inclusions with integer index  are those
arising from an outer action of a finite group $G$ on a
$II_1$ factor via fixed point algebras or crossed products. The index is $|G|$.
In the fixed point algebra case, $N=M^G\subset M$, it is well known that
$M_1$ can be identified with the crossed product
$M\rtimes G$ and   $N'\cap M_1$
with ${\mathbb C}G$ \cite{Jones1}.
This inclusion is irreducible and has depth $2$, the higher terms of the 
chain
$N'\cap N\subset N'\cap M\subset N'\cap M_1\subset\dots $
are determined  as in 4.7 a), \cite{GHJ}. It is well known that the
category ${\cal T}_M$ is described by the representation theory of $G$.
It is also well known that
this generalizes  to finite dimensional Kac algebras and that, by a result
 of Ocneanu, any irreducible finite index depth $2$ inclusion arises in
this way (see \cite{Szymanski}).
\medskip

\noindent b) {\it Bisch--Haagerup inclusions.}
  Recall that these
inclusions are obtained composing fixed point subfactors with crossed 
products subfactors:  $N=P^H\subset P\rtimes
K=M$ where $H$ and $K$ are finite groups acting properly and outerly on the
$II_1$ factor $P$. Recall from \cite{BH} the following results: $N\subset
M$ is always
extremal,
with integer index given by $|H||K|$,  the associated graph   is
amenable if and only if the group $G$ generated by $H$ and $K$ in 
$\text{Out}(P)$ is amenable, $N\subset M$ is
of finite depth if and only if $G$ is finite. Irreducible depth $2$ 
inclusions correspond
to matched pair of groups: $G=HK$, $H\cap K=\{e\}$. The corresponding Kac 
algebras
have been identified in \cite{HS}. Ocneanu's duality has been generalized 
by Nikshych and Vainerman to reducible depth $2$ inclusions. They proved 
that in this case
finite dimensional weak Hopf algebras (or quantum groupoids) in the sense 
of \cite{BS} replace Kac algebras
\cite{NV}. Vallin then proved that relative matched pairs (i.e. $G=HK$   
but $H\cap K$ is not required to be trivial) give rise to such inclusions, 
and hence to weak Hopf algebras
\cite{Vallin}.
\medskip

\noindent c) Known classes of extremal and {\it strongly amenable}
\cite{Popa} (hence
amenable) inclusions are given by: subfactors $N$ of the hyperfinite
$II_1$
factor $R$ with $[R:N]\leq4$ or of finite index and finite depth, e.g.
the
Jones subfactors \cite{GHJ} and 
Wenzl
subfactors arising from representations of the Hecke algebras of type $A$,
$B$, $C$, $D$ at roots of unity \cite{Wenzl1}, \cite{Wenzl2}. These 
provide non-embedable
categories as the index is not integral
\medskip

\noindent c) {\it Compact groups.}
Subfactors arising from
actions of
compact groups $G$ on $R$ and  finite dimensional unitary irreducible
$G$--representations have been considered in \cite{GHJ}, 
\cite{Wassermann88}. These, on the 
contrary, provide amenable embedable categories with integer values of the index.
\medskip

\noindent d) {\it Compact quantum groups.} 
Compact quantum groups give rise to examples of 
extremal non-amenable subfactors with 
embedable standard invariant, as shown by Banica in \cite{Banica2}. More 
precisely,   
he shows that, given a unitary f.d. representation $u$
of a compact quantum group, the 
selfintertwiners of the iterated tensor products  of $u$ with 
$\overline{u}$, gives rise to a standard 
$\lambda$--lattice in the sense of Popa, and hence, by the main result of
\cite{Popa2}, to an extremal 
inclusion of factors $N\subset M$ with $\lambda^{-1}=[M:N]=d(u)^2$,
where $d(u)$ is the quantum dimension of $u$. By construction, the corresponding lattice is embedded in the full lattice associated with the 
Hilbert space of $u$.  
However, by Theorem D of the same paper, this lattice is non-amenable 
if the quantum dimension of $u$ differs from its Hilbert space dimension, in agreement with Prop. 3.1.
\medskip

\noindent e) {\it Asaeda--Haagerup subfactors.} As is well known, the 
finite depth subfactors 
of indices $(5+\sqrt{13})/2$ and $(5+\sqrt{17})/2$ of \cite{AH} are not 
associated to 
classical groups or quantum groups. These provide 
examples of 
tensor categories ${\cal T}_M$ that can not be embedded into Hilbert 
spaces.
\medskip
\end{section}

\begin{section} {Statement of the results}

\noindent{\bf 4.1 Theorem} {\sl Let $N\subset M$ be a proper inclusion 
of
$II_1$
factors with finite index. For any integer $2\leq n\leq[M:N]$, let
$F\in M_n({\mathbb C})$ be an invertible matrix satisfying
$F\overline{F}=I,\quad \text{Trace}(F^*F)=[M:N].$
Then 
\begin{description}
\item{\rm a)}
there is an ergodic action of $A_o(F)$ on a unital
$C^*$--algebra ${\cal C}$ with spectral spaces
$\overline{L}_{u^{\otimes r}}=N'\cap M_{r-1}$, $r\geq0$, with inner 
product 
defined by the restriction of the normalized trace on $M_{r-1}$.
One has:
$$m(u^{\otimes r})=\text{dim}(N'\cap 
M_{r-1}),$$
where $u$ is the defining representation of $A_o(F)$.
\item{\rm b)}
In particular, if we choose $n=2$ we get an action of $S_{-\mu} U(2)$ with
$0<\mu\leq 1$ determined by $\mu+\mu^{-1}=[M:N]$.
\item{\rm c)}
If $[M:N]$ is an integer, 
say
$p$,
then we get an ergodic
$C^*$--action of the Kac compact quantum group $A_o(I_p)$. In this case
the unique invariant state is a trace.
\end{description}
}\medskip

\noindent{\it Remark} 
As is well known, a large class 
of subfactors have been constructed from compact groups \cite{GHJ}, 
\cite{Wassermann88}, quantum groups, first at roots of 
unity 
\cite{Wenzl1}, \cite{Wenzl2}, \cite{Wenzl3} and later even from compact 
quantum groups \cite{Banica2}. Now the modular theory of the Haar state of compact quantum 
groups is not trivial unless the coinverse is involutive \cite{Wcmp}.
Furthermore, in contrast to the classical results on group actions,
where an ergodic action implies a finite trace \cite{HLS},
or, for $SU(2)$, an algebra of type $I$  \cite{Wassermann1}, 
\cite{Wassermann2}, \cite{Wassermann3},
Wang showed that $A_u(Q)$ acts ergodically 
both on the type $III_\lambda$ factors of Powers, and on the Cuntz algebras, 
while $A_u(I_n)$ acts ergodically on the hyperfinite $II_1$ factor 
\cite{Wang2}. Thus we cannot expect our $C^*$--algebra, ${\cal C}$, to be finite, 
in general, and we do not know how the type of ${\cal C}$ is related to 
properties of the spectral functor of the action.
\medskip

The following result is related to Prop. 3.1. For simplicity, we give a direct proof.\medskip

\noindent{\bf 4.2  Corollary} {\sl If the inclusion $N\subset M$ is 
extremal and 
amenable in the sense of  \cite{Popa} and if $[M:N]$ is not an integer then the                above ergodic action of
 $A_o(F)$ as is not embedable into the translation action on $A_o(F)$.}\medskip

\noindent{\it Proof}  By  \cite{Popa},
$$\lim \text{dim}(\overline{L_{u^{\otimes r}}})^{1/r}=\lim\text{dim}(N'\cap M_{r-1})^{1/r}=[M:N].\eqno(4.1)$$
Hence if $[M:N]$ is not 
an integer, $\text{dim}(\overline{L_{u^{\otimes 
r}}})>n^r$ for $r$ large enough. Thus  the action cannot be embedable in 
the translation action, as this would imply
$\overline{L_{u^{\otimes r}}}\subset H_{u^{\otimes r}}$ for all $r$.
\medskip

 Let us discuss some examples.
\medskip

\noindent {\it Irreducible depth $2$ inclusions.} We compare our
construction with Ocneanu's duality recalled in a), which reconstructs an 
outer action of a
f.d. Kac algebra $G$ from a depth $2$ irreducible finite index
inclusion.  The  fixed point subfactor 
gives rise to a category ${\cal T}_M$ isomorphic to the category generated 
by the
tensor powers of the regular representation $\lambda$ of $G$.
Since
$[M:M^G]=\text{dim}(\ell^2(G))=:n$, $F:=I_n$, and
hence $A_o(n)$, is a natural
choice.
Hence the spectral space $\overline{L}_{u^{\otimes r}}$ 
 of the resulting ergodic action of $A_o(n)$,
 is given by the space of
 fixed vectors of $\ell^2(G)^{\otimes r}$ under
the $r$-th tensor power of $\lambda$.
On the other hand the regular representation $\lambda$ of $G$ is a real
object
of intrinsic dimension $n$. Hence $G$,
regularly represented, may be
regarded as a quantum subgroup of $A_o(n)$. Our construction
thus gives the quantum quotient space $G\backslash A_o(n)$. \medskip

\noindent{\it The case $[M:N]=2$.}
 In this case $N=M^{{\mathbb Z}_2}$ by Goldman's theorem \cite{Goldman}. 
 $A_o(2)$ is the only possible quantum group arising in our 
framework, up to similarity between compact matrix quantum groups. 
In particular, if  $u=(u_{ij})$ and $v=(v_{ij})$ denote the fundamental 
representations of 
 $A_o(2)$ and $S_{-1}U(2)$ respectively, the map $\phi:C(S_{-1}U(2))\to 
A_o(2)$ 
defined by $\iota\otimes \phi(v)=V uV^{-1}$, 
with $V=(V_{ij})$ the scalar valued matrix  $V_{11}=-V_{22}=i$, 
$V_{12}=-V_{21}=1$, is a  natural similarity \cite{Banica}. 
 We next identify the corresponding quotient space ${\mathbb 
Z}_2\backslash S_{-1}U(2)$.
It is clear from the work of \cite{Tomatsu}, based on previous results of 
\cite{Podles}, that ${\mathbb Z}_2$ gives rise to two 
different quotient $S_{-1}U(2)$--spaces, one corresponding to the usual
diagonal embedding of ${\mathbb Z}_2$, and another one.
We show that our ergodic action identifies with the latter.
By an argument of \cite{Tomatsu}, it suffices to show that the 
corresponding restriction 
map $r: C(S_{-1}U(2))\to A_o(2)\to C({\mathbb Z}_2)$
satisfies
$r(v_{11})(g)=0$.
We have, up to a scalar,
 $\phi(v_{11})=u_{11}-iu_{21}+iu_{12}+u_{22}$.  The restriction of $u$
to  the subgroup  $\lambda({\mathbb Z}_2)$ is its fundamental representation 
$\hat{u}=(\hat{u}_{ij})$ as a subgroup of $U(\ell^2({\mathbb Z}_2))$, which is given by convolution   on $\ell^2({\mathbb Z}_2)$. Hence
 $\hat{u}_{11}(g)=\hat{u}_{22}(g)=0$, $\hat{u}_{21}(g)=\hat{u}_{12}(g)=1$,
and the proof is complete.
\medskip

\noindent{\it Bisch--Haagerup subfactors for relative matched pairs.}
Let us consider two finite subgroups $H$, $K$
of $\text{Out}(R)$ forming a relative matched pair,
as recalled in b) of the previous section.
Our construction realizes Vallin's  quantum groupoid \cite{Vallin}
associated to $R^H\subset R\rtimes K$
as a virtual subgroup
of $A_o(n)$, with $n=|H||K|$.\medskip

\noindent {\it Jones subfactor.} Consider
the Jones subfactor $R_\beta\subset R$ of the hyperfinite $II_1$ 
factor with index $\beta=4\cos^2\pi/m$, with $m \geq 4$. Here  
${R_\beta}'\cap R_{r-1}$ is the algebra $B_{\beta,r}$
 generated by the Jones projections 
$e_1,\dots, e_{r-1}$.
It carries  a unitary representation of the braid group  ${\mathbb B}_r$ as we are in the 
case $\beta<4$,
see, \cite{Jones1} or  \cite{GHJ}
(hence 
the tensor $C^*$--category 
generated by the Hilbert bimodule ${}_N M_N$ has a unitary braiding, as
 we find such representations in the 
$C^*$--algebras $(M^{\otimes r},M^{\otimes r})\simeq {R_\beta}'\cap 
R_{2r-1}$).
If we apply the previous theorem, we deduce that $B_{\beta,r}$, 
regarded as a Hilbert space with inner product defined by its Markov 
trace, does arise as the spectral space $\overline{L_{u^{\otimes r}}}$
of an ergodic $C^*$--action of  
$S_{-\mu} U(2)$, with
$\mu+\mu^{-1}=\beta$. 
\medskip

\noindent {\it Wenzl subfactors.} In more generality, let $N\subset M$ be  the 
$II_1$ subfactor 
arising from 
quantum groups at roots of unity as in \cite{Wenzl3}. The higher 
relative commutants $N'\cap M_{r-1}$ are there shown to correspond to 
the arrow spaces of the fusion tensor category
of the quantum group  (Theorem 4.4 in 
\cite{Wenzl3}). These spaces, by the previous theorem, again arise as spectral 
spaces 
of ergodic actions of compact quantum groups.
Furthermore the  quantum dimension of the compact quantum group in question depends on the roots of unity. Thus these algebraic quantum groups seem to be virtual
quantum subgroups of compact quantum groups.

\medskip

\noindent {\it Banica  subfactors.} 
Banica's subfactors associated to a unitary representation $v$ of a 
compact 
quantum group have relative commutant $N'\cap M_r$ given by the 
selfintertwiners of the tensor product representation 
$v\otimes\overline{v}\otimes v\dots$ ($r+1$ factors) \cite{Banica2}. By Frobenius reciprocity, this space is linearly isomorphic to the space of invariant vectors
of the representation $(v\otimes\overline{v})^{\otimes r+1}$.
Hence 
this space
may be regarded as the spectral space $\overline{L}_{u^{\otimes r+1}}$ of an ergodic action of $A_o(F)$.\medskip

For convenience, we next give a presentation by generators and relations of
 the ergodic $C^*$--algebra ${\cal C}$  in terms of the higher relative commutants.  
 This will be deduced from a simpler presentation in terms of natural relations in Ocneanu's tensor $C^*$--category of Hilbert bimodules that will appear  in the course of the proof in the next section, (cf. relations  $(5.1)$--$(5.3)$ and Theorem 5.2).
 
 Before stating the result we need some notation.
Set $H:={\mathbb C}^{\times n}$, $j:=Fc$, where  $c:H\to H$ is  the antiunitary 
fixing the canonical basis of $H$, and 
$$R_u:=\sum\psi_k\otimes j\psi_k,$$
where  $(\psi_k)$ is  any orthonormal basis of $H$. We  
introduce 
certain reduced words 
in the algebra generated by the Jones projections.
Set $\lambda:=[M:N]^{1/2}$
and for nonnegative integers $k$, $r$, $s$, define elements
$p^{(k)}_{r,s}\in M_{k+r+s-1}$ by
$$p^{(k)}_{0,s}=p^{(k)}_{r,0}:=I,$$
and, for $r,s\geq1$,
$$p^{(k)}_{r,s}:=\lambda^{rs}(e_{r+k}e_{r+k-1}\dots 
e_{1+k})\dots(e_{r+k+s-1}e_{r+k+s-2}\dots e_{s+k}).\eqno(4.2)$$
 We shall 
simply write $p_{r,s}$
for $p^{(0)}_{r,s}$. By \cite{Pimsner-Popa},  $p_{r,r}$ reduces to a scalar multiple of the Jones projection associated to   $N\subset M_{r-1}$.

In the next result, in order to avoid confusion with the 
tensor products, we shall denote 
 the $r$-th tensor 
powers of $H$ and $u$
 by $H^r$ and $u^r$ respectively.
\medskip

\noindent{\bf 4.3 Theorem} {\sl
The $C^*$--algebra ${\cal C}$ is obtained by completing the $^*$--subalgebra with generators
$\overline{T}\otimes\xi$, $T\in N'\cap M_{r-1}$, $\xi\in H^r$,
$r=0,1,2,\dots$,
and relations, where $T'\in N'\cap M_{s-1}$, $\xi'\in H^s$, 
$\xi_1,\dots,\xi_r\in H$, $\eta\in H^{r+s+2}$, $\eta'\in H^{r+s}$,
\begin{description}
\item{\rm a)}
$(\overline{T}\otimes\xi)(\overline{T'}\otimes\xi')=\overline{Tp_{r,s}T'}
\otimes \xi\xi',$
\item{\rm b)}
$(\overline{T}\otimes\xi_1\dots\xi_r)^*=\overline{T^*}\otimes
j\xi_r\dots j\xi_1,$
\end{description}
for $r\geq s$:
\begin{description}
\item{\rm c)}
$\overline{S}\otimes
(1_{u^r}\otimes R_u^*\otimes 1_{u^s}\eta)=
\lambda \overline{Sp^{(2s)}_{r-s,2}}\otimes\eta$,
\item{\rm c')}
$\overline{S'}\otimes (1_{u^r}\otimes
R_u\otimes 1_{u^s}\eta')=
\lambda\overline{E_{r+s}E_{r+s+1}(S'(p^{(2s)}_{r-s,2})^*)}\otimes\eta'$,
\end{description}
for $r<s$:
\begin{description}
\item{\rm d)}
$\overline{S}\otimes(1_{u^r}\otimes R_u^*\otimes 1_{u^s}\eta)=
\lambda \overline{p^{(2r)}_{2,s-r}S}\otimes\eta$,
\item{\rm d')}
$\overline{S'}\otimes (1_{u^r}\otimes
R_u\otimes 1_{u^s}\eta')=
\lambda
\overline{E_{r+s}E_{r+s+1}((p^{(2r)}_{2,s-r})^*S')}\otimes\eta'$,
\end{description}  
 in the maximal
$C^*$--norm.
The $A_o(F)$--action $\beta$ is uniquely defined by
$$\beta(\overline{T}\otimes\xi)=\overline{T}\otimes u^{\otimes r}(\xi),\quad  T\in 
N'\cap M_{r-1},\xi\in H^r,$$
where $u$ is the defining representation of $A_o(F)$ on $H$.
}\medskip

\end{section}

\begin{section}{Proof of the results}

We shall refer to \cite{DRInventiones} for the definition of an abstract 
tensor 
$C^*$--category ${\cal T}$. 
The tensor unit object will be denoted $\iota$.
We shall always assume $(\iota,\iota)={\mathbb C}$.
An object $\rho$ will be called real (or pseudoreal) if there is an
$R\in(\iota,\rho^2)$ satisfying $R^*\otimes 1_\rho\circ 1_\rho\otimes R=1_\rho$ (or $R^*\otimes 1_\rho\circ 1_\rho\otimes R=-1_\rho$). 
\medskip

\noindent{\bf 5.1 Theorem} {\sl Let $\rho$ be a real or pseudoreal object of ${\cal T}$  defined by
$R\in(\iota,\rho^2)$ with $\|R\|^2\geq2$. For any integer
$2\leq n\leq\|R\|^2$ let $F\in M_n({\mathbb C})$ be an invertible matrix such that $F\overline{F}=\pm I$ and $\text{Trace}(F^*F)=\|R\|^2$. Then there is an ergodic action of $A_o(F)$ on a unital $C^*$--algebra 
${\cal C}$ with spectral functor $\overline{L}$, where  
$\overline{L}_{u^{\otimes r}}=(\iota,\rho^r)$ and 
$\overline{L}_{\sum_k\psi_k\otimes F\psi_k}$ is left compositon by $R$. In 
particular, for $n=2$ 
we get an action of
 $S_{\mu} U(2)$ for a nonzero
$-1\leq\mu\leq 1$ determined by $|\mu+\mu^{-1}|=\|R\|^2$ and $\mu>0$ if   $\rho$ is pseudoreal and negative 
otherwise.}\medskip

The above theorem was proved in \cite{PR}. We shall outline the proof for convenience. 
\medskip

\noindent{\it Outline of proof.}
One first shows that the tensor $^*$--subcategory 
 generated by an arrow $R$ making $\rho$ real 
(or pseudoreal)  is uniquely determined by the quantum dimension $d=\|R\|^2$.
If  $F$ is  chosen as indicated, it gives rise to a realization  of this category in the category of Hilbert spaces which generates the representation category of $A_o(F)$. 
Thus we have a tensor 
$^*$--functor, still denoted by $\rho$,  from the full tensor subcategory 
of $\text{Rep}(A_o(F))$ generated by
the fundamental representation $u$ 
 to ${\cal T}_\rho$ taking $u$ to 
$\rho$ 
and the basic intertwiner $\sum_k\psi_k\otimes F\psi_k$ to $R$. We next apply 
 the duality theorem in \cite{PR} (cf. Sect. 2) to the quasitensor functor
$u^{\otimes r}\to (\iota, \rho^r)$. That theorem yields
 a presentation by generators and relations of  the dense linear space (in fact a $^*$--subalgebra) ${\cal C}_{\text{sp}}$ generated
 by the spectral elements  $T^*\otimes\xi$ with relations: 

$$({T^*}\otimes\xi)({T'}^*\otimes\xi')={T^*\otimes{T'}^*}
\otimes\xi\xi',\eqno(5.1)$$
$$({T^*}\otimes \xi)^*=\rho(R_{u^r}^*)\circ1_{\rho_{\overline{u^r}}}\otimes 
T\otimes j_{u^r}\xi,\eqno(5.2)$$
$$T^*{\rho(A)}\otimes\xi'={T^*}\otimes 
A\xi',\eqno(5.3)$$
where  $T\in(\iota,\rho^r)$,  $T'\in(\iota,\rho^s)$, $\xi\in H^r$, $\xi'\in H^s$, $A\in(u^s, u^r)$,  
and $j_{u^r}$ 
defines a solution $R_{u^r}$, $\overline{R}_{u^r}$ of the conjugate 
equations for the $r$-th 
tensor power $u^r$ of the defining representation  $u$.
\medskip

Finiteness of the Jones index implies that ${\cal T}_M$ has conjugates.
We start by recalling some facts on module bases for $M_N$ from 
\cite{Pi-Po}.
A  basis for $M$, as a right Hilbert module, is a finite 
set of elements
$(u_i)$ in $M$ such that $\sum_i u_iE(u_i^*m)=m$ for all $m\in M$.
We shall 
refer to $(u_i)$ as a Pimsner-Popa basis.
Such bases exist as $M_N$ is finitely generated.
If $(u_i)$ is a Pimsner-Popa basis then $\sum_i 
u_ie_1u_i^*=I$ in $M_1$. We have: $\sum_i u_iu_i^*=[M:N].$ The following 
result has been shown in \cite{KPW} in more generality. 
\medskip

\noindent{\bf 5.2 Theorem} {\sl 
${}_NM_N$ is a real object of the category ${\cal T}_M$.  A 
solution of the conjugate 
equations for $M$ is given by
$$R=\overline{R}=\sum_i u_i\otimes u_i^*,\eqno(5.4)$$
where $(u_i)$ is a Pimsner-Popa basis for $M$.
One has $\|R\|^2=[M:N]$.
}\medskip

\noindent{\it Proof} As we have a tensor product over $N$
it is easily checked that $R$ is independent of the choice of the Pimsner--Popa basis. 
On the other hand for any unitary $u\in N$, $(uu_i)$ is another 
Pimsner--Popa basis, hence $uRu^*=R$, showing that $R$ is an intertwiner 
in ${\cal T}_M$. For $m\in M$,
$$R^*\otimes 1_M\circ1_M\otimes R(m)=\sum_iR^*\otimes  
1_M(m\otimes u_i\otimes u_i^*)=$$
$$\sum_{i,j}<u_j\otimes u_j^*, m\otimes u_i> u_i^*=\sum_{i,j}<u_j^*, 
E(u_j^*m)u_i>u_i^*=$$
$$\sum_{i,j}E(u_jE(u_j^*m)u_i)u_i^*=\sum_i E(mu_i)u_i^*=m.$$
We also have
$$\|R\|^2=R^*R=\sum_{i,j}<u_i\otimes u_i^*, u_j\otimes u_j^*>=$$
$$\sum_{i,j}E(u_iE(u_i^*u_j)u_j^*)=\sum_j E(u_ju_j^*)=[M:N].$$
\medskip

Applying Theorem 5.1 to  Ocneanu's 
category ${\cal T}_M$
and the intertwiner $R=\sum u_i\otimes u_i^*$,  we get the desired ergodic action of $A_o(F)$ 
on a unital
$C^*$--algebra with spectral spaces $\overline{L}_{u^{\otimes 
r}}=(\iota,M^{\otimes r})\simeq N'\cap M_{r-1}$.
It remains to identify
the inner products on the spaces $(\iota,M^{\otimes r})$ 
arising from the category ${\cal T}_M$ and the algebraic presentation of 
${\cal C}_{\text{sp}}$ in terms of the higher relative 
commutants $N'\cap M_{r-1}$. 

Now the inner product of $(\iota,M^{\otimes r})$ is the restriction of the 
$N$--valued inner product of 
$M^{\otimes r}$ which, through  the unitary  
identification $M^{\otimes r}\simeq L^2(M_{r-1})$ as 
$N$--$N$--correspondences in the sense of Connes \cite{Connes},  arises from the 
normalized trace 
of 
$M_{r-1}$ (see, e.g.,  Prop.\ 3.1 in  \cite{Bisch}).

For later use, we shall need  explicit Hilbert $N$--bimodule unitaries from 
$M^{\otimes r}$ to $M_{r-1}$. 
 Consider the $N$--bimodule isomorphism between
$M_{r+1}:=M_re_{r+1}M_r$ and
$M_r\otimes_{M_{r-1}}M_r$,
given by
$$m\otimes m'\to\lambda me_{r+1} m',\quad m,m'\in M_r.$$
Regard 
$M_r$ as a $M_{r-1}$--Hilbert bimodule with inner product defined by the 
normalized
conditional expectation $E_r: M_r\to M_{r-1}$.
Then the tensor product
of the $M_{r-1}$--valued inner products on
$M_r\otimes_{M_{r-1}}M_r$,
corresponds, under the above
isomorphism $U$
to the inner product induced by the (normalized) conditional expectation
$E_rE_{r+1}: M_{r+1}\to M_{r-1}$
by
$$<S, T>=E_rE_{r+1}(S^*T), \quad S, T\in M_{r+1}.$$
Iterating this procedure (recall  $(3.5)$), leads to
  canonical isomorphisms of $N$--bimodules
 $\Gamma_r: M^{\otimes r}\to M_{r-1}$
 transforming
the tensor product of $N$--valued inner products into the inner product
induced by the conditional expectation
$$E_{(r)}:=EE_1\dots E_{r-2}E_{r-1}:
M_{r-1}\to N.\eqno(5.5)$$
 $\Gamma_r:M^{\otimes r}\simeq M_{r-1}$ is obtained in the following way. 
First   
replace  each factor $M$ occurring in $2$nd to $r-1$th position in
$M^{\otimes r}$ by $M\otimes_M M$,
giving  $2r-2$ factors tensored alternately over $N$ and
$M$.
 Thus $$M^{\otimes r}=(M\otimes_N M)^{\otimes_M r-1}.$$
Finally use the isomorphism $M\otimes_N M\simeq M_1$ giving by iteration
$$M^{\otimes r}\simeq M_1^{\otimes_M r-1}\simeq M_2^{\otimes_{M_1}
r-2}\simeq\dots\simeq M_{r-2}\otimes_{M_{r-3}} M_{r-2}\simeq M_{r-1}.$$
The resulting isomorphism  $\Gamma_r $   is
$$m_1\otimes\dots\otimes
m_r\to\lambda^{<r>}
m_1e_1m_2e_2e_1m_3\dots m_{r-1}e_{r-1}\dots
e_1m_r,\eqno(5.6)$$
where $<r>:=(r-1)+(r-2)+\dots+1=r(r-1)/2$.
Summarizing,  one has the following result.
\medskip

\noindent{\bf 5.3 Proposition} {\sl  
Under the $N$--bimodule isomorphisms $\Gamma_r: M^{\otimes r}\to 
M_{r-1}$, the $N$--valued tensor power inner product
on $M^{\otimes r}$ corresponds to the inner product on $M_{r-1}$
defined by the trace-preserving conditional expectation $E_{(r)}: M_{r-1}\to N$:
$$<S,T>:=E_{(r)}(S^*T),\quad S,T\in M_{r-1}.$$
In particular, the Hilbert space
structure of $(\iota,M^{\otimes r})$
defined by the category ${\cal T}_M$  
corresponds, under the restriction of
$\Gamma_r$,
to the inner product of  $N'\cap
 M_{r-1}$ defined by the restriction of the normalized trace on $M_{r-1}$.
}\medskip

We need to establish the algebraic presentation of the dense spectral 
$^*$--subalgebra ${\cal C}_{\text{sp}}$. From   \cite{PR} (cf.\ Theorem 5.1 and its proof)
applied 
to    ${\cal T}={\cal T}_M$,  Ocneanu's tensor 
$C^*$--category 
associated with the inclusion $N\subset M$ as in section 3, and to the 
tensor 
$^*$--functor 
from the full tensor $C^*$--subcategory of $\text{Rep}(A_o(F))$
generated by the defining representation $u$ to ${\cal T}_M$ and taking 
$u$ to ${}_NM_N$
and the basic intertwiner $\sum_k\psi_k\otimes F\psi_k$ to $R=\sum u_i\otimes u_i^*$,
we 
know that an algebraic 
presentation of ${\cal C}_{\text{sp}}$ is given by generators 
${T^*}\otimes\xi$
with $T\in(\iota, M^{\otimes r})$, $\xi\in H^r$ subject to the relations 
$(5.1)$--$(5.3)$.
 
 We start by computing the $^*$--involution.
Starting from the antilinear intertwiner $j_u:=j=Fc$, corresponding to
$R_u=\sum\psi_i\otimes j\psi_i$, for $r=1$, we can 
form its tensor power $j_{u^r}(\xi_1\dots\xi_r):=j(\xi_r)\dots j(\xi_1)$. With this 
choice we have 
$$R_{u^r}=1_{u^{r-1}}\otimes R_u\otimes 1_{u^{r-1}}\circ 
R_{u^{r-1}}.$$
By tensoriality of 
the inclusion $\rho$ of the full tensor 
$C^*$--subcategory of
$\text{Rep}(A_o(F))$ in 
 $ {\cal 
T}_M$,  
$$\rho(R_{u^r})=1_{M^{\otimes r-1}}\otimes \rho(R_u)\otimes 
1_{M^{\otimes r-1}}\circ 
\rho(R_{u^{r-1}})=$$
$$\sum_{i_1,\dots i_r}u_{i_1}\otimes\dots 
u_{i_r}\otimes{u_{i_r}}^*\otimes\dots{u_{i_1}}^*\in M^{\otimes 2r}.$$
On the other hand, since $u^r$ is real, it is selfconjugate, hence 
$\overline{u^r}=u^r$. Therefore, for $T\in(\iota, M^{\otimes r})$,
$$(\overline{T}\otimes\xi_1\dots\xi_r)^*=\overline{1_{M^{\otimes 
r}}\otimes 
T^*\circ\rho(R_{u^r})}\otimes j\xi_r\dots j\xi_1=$$
$$\overline{u_{i_1}\otimes\dots\otimes u_{i_r}<T, 
u_{i_r}^*\otimes\dots\otimes u_{i_1}^*>}\otimes j\xi_r\dots j\xi_1.$$
We thus need to compute $\sum u_{i_1}\otimes\dots\otimes
u_{i_r}<T,u_{i_r}^*\otimes\dots\otimes u_{i_1}^*>$ for
$T\in(\iota, M^{\otimes r})$ under the identification 
$(\iota, M^{\otimes r})\simeq N'\cap
M_{r-1}$.
Now
it is easy
to see, using the trace norm of $N$, that $E(SX)=E(XS)$ for $S\in
N'\cap M$ and $X\in M$, hence,  for $r=1$, and $T\in N'\cap M$,
$$\sum_i u_i<T, u_i^*>=\sum_i u_iE(T^*u_i^*)=
\sum_i u_iE({u_i}^*T^*)=T^*.$$
In the general case, if $T\in N'\cap M_{r-1}$, identifying  $M^{\otimes r}$  and $M_{r-1}$ as right Hilbert
$N$--modules, $u_{i_1}\otimes\dots\otimes u_{i_r}$ corresponds to a
Pimsner--Popa basis for the inclusion $N\subset M_{r-1}$, and the above
argument for $r=1$ shows that
$\sum u_{i_1}\otimes\dots\otimes u_{i_r}<T, u_{i_r}^*\otimes\dots\otimes
u_{i_1}^*>$ corresponds to the adjoint $T^*$ of $T$ in $N'\cap M_{r-1}$.
Summarizing:   
\medskip

\noindent{\bf 5.4 Proposition} {\sl Under the restriction of the 
$N$--bimodule 
unitary
$\Gamma_r: M^{\otimes r}\to M_{r-1}$, the antilinear map
$$J_r:T\in(\iota,M^{\otimes r})\to 1_{\rho_{\overline{u^{\otimes 
r}}}}\otimes
T^*\circ\rho(R_{u^{\otimes r}})\in(\iota, M^{\otimes r})\eqno(5.7)$$ 
corresponds 
to the antiunitary
$^*$--involution
$$T\in N'\cap M_{r-1}\to T^*\in N'\cap M_{r-1}.\eqno(5.8)$$
Therefore for the resulting ergodic $C^*$--action of 
$A_o(F)$ 
we have:
$$m(u^r)=\text{dim}(N'\cap M_{r-1}),\eqno(5.9)$$
where  $m$ is the quantum multiplicity.
}\medskip

Hence formula b)  in Theorem 4.3  has been established.

As the quantum multiplicity necessarily takes on its minimal value, 
$II_1$ subfactors do not yield all canonical ergodic actions of $A_o(F)$.

Before  exploiting the operator product $(5.5)$ in ${\cal C}$, we complete the proof of Theorem 4.1.\medskip

\noindent{\bf Proof of Theorem 4.1 c)}
\medskip

We need to show that the unique invariant state 
$h$   is a trace. In the special case that we are 
considering, $A_o(I_p)$ has an involutive coinverse \cite{WVD}, hence $j_w$ 
can be chosen  antiunitary for all representations $w$.
In particular the  
solution $R_{u^r}$, $\overline{R}_{u^r}$ of the conjugate equations is 
standard. By  Cor. A.10 of \cite{PR2}, the  solution $\hat{R}_{u^r}$, 
$\hat{\overline{R}}_{u^r}$ is also standard. 
 We close ${\cal T}_M$ under subobjects and then
extend the inclusion functor $\rho$ to a relaxed tensor $^*$--functor from $\text{Rep}(A_o(F))$ to this closure, so that 
$\rho_v$ is now defined for any   irreducible 
subrepresentation $v$ of some $u^r$. It follows  that  $J_r\upharpoonright_{(\iota,\rho_v)}$
coincides with an antilinear invertible $J_v$ constructed from  a  
normalized
solution
of the conjugate equations for $v$. Since $J_r$ is antiunitary, $J_v$ is 
antiunitary as well.  Let $\hat{\rho}$ be a complete set of irreducible representations of $A_o(F)$.  If $v,w\in\hat{\rho}$, $a=\overline{S}\otimes\psi$, 
$b=\overline{T}\otimes\phi$, $\psi\in H_v$, $\phi\in H_w$, 
$S\in(\iota,\rho_v)$, $T\in(\iota,\rho_w)$ then $v\otimes w$ contains the 
trivial representation $\iota$ if and only if $w=\overline{v}$, and the 
multiplicity of $\iota$ is $1$. Hence 
$S=\|\overline{R_v}\|^{-1}\overline{R_v}\in(\iota, 
v\otimes\overline{v})$ is an isometry.
It follows that
$$h(ab)=
\delta_{w,\overline{v}}\|\overline{R_v}\|^{-2}
\overline{\rho(\overline{R_v}^*)S\otimes 
T}\otimes\overline{R_v}^*\psi\otimes \phi=$$
$$\delta_{w,\overline{v}}\|\overline{R_v}\|^{-2}
<J_vS, T><j_v\psi,\phi>.$$
Similarly,
$$h(ba)=
\delta_{v,\overline{w}}\|\overline{R_w}\|^{-2}
<J_wT, S><j_w\phi,\psi>$$
which in turn equals
$$\delta_{w,\overline{v}}\|{R_v}\|^{-2}
<J_v^{-1}T, S><j_v^{-1}\phi,\psi>=$$
$$\delta_{w,\overline{v}}\|{R_v}\|^{-2}
<J_vS, T><j_v\psi,\phi>=h(ab).$$
\medskip

We next  
spell out the tensor product operation between arrows
in Ocneanu's category in terms of the higher relative commutants.

In detail, the isomorphism $\Gamma_r$ allows us to write down the 
tensor product
$\xi\otimes\eta$ of elements $\xi\in M^{\otimes r}$ and
$\eta\in M^{\otimes s}$ in terms of a  bilinear map
$$M_{r-1}\times M_{s-1}\to M_{r+s-1}.$$
If, e.g.,  $\xi=\xi_1\otimes \xi_2$,
$\eta=\eta_1\otimes \eta_2\otimes \eta_3$,  
$$\Gamma_5(\xi\otimes\eta)=\lambda^{10}\xi_1
e_1\xi_2
e_2e_1\eta_1e_3e_2e_1\eta_2e_4e_3e_2e_1\eta_3=$$
$$\lambda^{9}\Gamma_2(\xi)e_2e_1\eta_1e_3e_2e_1\eta_2e_4e_3e_2e_1\eta_3=$$
$$\lambda^9\Gamma_2(\xi)
e_2e_1e_3e_2e_4e_3
(\eta_1e_1\eta_2e_2e_1\eta_3)=$$
$$\lambda^{6}\Gamma_2(\xi)e_2e_1e_3e_2e_4e_3\Gamma_3(\eta).$$   
In general, if $\xi\in M^{\otimes r}$, $\eta\in M^{\otimes s}$
then
$$\Gamma_{r+s}(\xi\otimes\eta)=\lambda^{<r+s>-<r>-<s>}\Gamma_r(\xi)
(e_r\dots e_1)(e_{r+1}\dots e_2)\dots(e_{r+s-1}\dots e_s)
\Gamma_s(\eta)=$$
$$\lambda^{rs}\Gamma_r(\xi)
(e_r\dots e_1)(e_{r+1}\dots e_2)\dots(e_{r+s-1}\dots 
e_s)\Gamma_s(\eta)=\Gamma_r(\xi)p_{r,s}\Gamma_s(\eta),$$
with $p_{r,s}=p^{(0)}_{r,s}$.
\medskip

\noindent{\bf 5.5 Proposition} {\sl
Under the $N$--bimodule unitary $\Gamma_r: M^{\otimes  r}\to M_{r-1}$
the tensor product
$$\xi\in(\iota, 
M^{\otimes r}),\eta\in(\iota, M^{\otimes s})\to\xi\otimes\eta\in(\iota,
M^{\otimes r+s})\eqno(5.10)$$
defined in ${\cal T}_M$ corresponds to the map
$$S\in N'\cap M_{r-1}, T\in N'\cap M_{s-1}\to Sp_{r,s}T\in
N'\cap M_{r+s-1}.\eqno(5.11)$$}\medskip

We are left to account for the $\text{Rep}(A_o(F))$--bimodule structure 
$(5.3)$ used to define ${\cal C}_\rho$. 
Now the  full tensor $^*$--subcategory of $\text{Rep}(A_o(F))$ generated by $u$  is generated, as a
linear category, by the arrows $1_{u^{\otimes r}}\otimes R_u\otimes 
1_{u^{\otimes
s}}$, $r,s\geq0$, and by their adjoints. Therefore we are led to compute
the composition of arrows of the form
$1_{M^{\otimes r}}\otimes R\otimes 1_{M^{\otimes s}}\circ T$
and
$1_{M^{\otimes r}}\otimes R^*\otimes 1_{M^{\otimes s}}\circ S$
with $T\in(\iota, M^{\otimes r+s})$, $S\in(\iota, M^{\otimes r+s+2})$, 
$r,s\geq0$,
in terms of the higher relative commutants.
For $\xi\in M^{\otimes  r}$, $\eta\in M^{\otimes s}$,
$$1_{M^{\otimes r}}\otimes R\otimes1_{M^{\otimes 
s}}\xi\otimes\eta=\xi\otimes
R\otimes\eta,$$
hence
$$\Gamma_{r+2+s}(1_{M^{\otimes r}}\otimes
R\otimes1_{M^{\otimes s}}(\xi\otimes\eta))=$$
$$\Gamma_{r+2}(\xi\otimes
R)p_{r+2,s}\Gamma_s(\eta)=$$
$$\Gamma_r(\xi)p_{r,2}\Gamma_2(R)p_{r+2,s}\Gamma_s(\eta).$$
Now
$$\Gamma_2(R)=\sum_i\Gamma_2(u_i\otimes u_i^*)=\lambda\sum_i
u_ie_1u_i^*=\lambda.$$
Therefore
$$\Gamma_{r+2+s}(1_{M^{\otimes r}}\otimes
R\otimes1_{M^{\otimes s}}(\xi\otimes\eta))=
\lambda\Gamma_r(\xi)p_{r,2}p_{r+2,s}\Gamma_s(\eta).\eqno(5.12)$$
We therefore need to write $p_{r,2}p_{r+2,s}$ as a reduced word in the  
algebra generated by the Jones projections $e_1,\dots, e_{r+s+1}$.
It is known that elements of the form $$(e_{j_1}e_{j_1-1}\dots
e_{i_1})(e_{j_2}\dots e_{j_2-1}\dots e_{i_2})\dots(e_{j_p}e_{j_p-1}\dots
e_{i_p})$$
are in reduced form if $j_1<j_2<\dots<j_p$ and $i_1<i_2<\dots<i_p$
\cite{GHJ}. The following lemma is useful.
\medskip

\noindent{\bf  5.6  Lemma} {\sl For 
$p\leq j\leq r<s$ we 
have
$$(e_{r}e_{r-1}\dots e_j)(e_se_{s-1}\dots e_p)=
\lambda^{-2}(e_r\dots
e_p)(e_s\dots e_{j+2}) \quad\text{for } s>j+1,$$
$$(e_{r}e_{r-1}\dots e_j)(e_se_{s-1}\dots e_p)=
\lambda^{-2}(e_r\dots
e_p) \quad\text{for } s=j+1.$$
}\medskip

\noindent{\it Proof} We do the computation only in the case $s>j+1$.
$$(e_{r}e_{r-1}\dots e_j)(e_s\dots e_{j+2}e_{j+1}e_j\dots e_p)=
(e_{r}e_{r-1}\dots e_{j+1})(e_s\dots e_{j+2}e_je_{j+1}e_j\dots e_p)=$$
$$\lambda^{-2}(e_{r}e_{r-1}\dots e_{j+1})(e_s\dots e_{j+2}e_j\dots e_p)=
\lambda^{-2}(e_{r}e_{r-1}\dots e_{j+1}e_j\dots e_p)(e_s\dots e_{j+2}).$$
\medskip

\noindent{\bf 5.7 Lemma} {\sl For $s\geq0$ and $r\geq s$ we have:
$$p_{r,2}p_{r+2,s}=p_{r,s}p^{(2s)}_{r-s,2}.$$}\medskip

\noindent{\it Proof} The formula is obvious for $s=0$. We can then assume
$s>0$. We write down  the left hand side explicitly:
$$p_{r,2}p_{r+2,s}=\lambda^{2r+2s+rs}
(e_r\dots e_1)
(e_{r+1}\dots e_2)(e_{r+2}\dots e_1)\dots(e_{r+s+1}\dots e_s).$$
We have $2+s$ factors between parentheses. Let us  apply the
previous lemma
iteratively between the second and the third factor.
If $r+2>3$, i.e. $r>1$, we have, 
$$(e_{r+1}\dots e_2)(e_{r+2}\dots e_1)=\lambda^{-2}(e_{r+1}\dots
e_1)
(e_{r+2}\dots e_{2+2}).$$
We proceed to apply the lemma to the new third and old fourth
factors: if $r+3>5$, i.e. $r>2$,
$$(e_{r+2}\dots e_{2+2})(e_{r+3}\dots e_2)=\lambda^{-2}(e_{r+2}\dots
e_2)(e_{r+3}\dots e_{2+4}).$$ 
If $n>s$
after $s$ iterations of the lemma we still find
a product of $2+s$ factors:
$$p_{r,2}p_{r+2,s}=\lambda^{2r+rs}(e_r\dots e_1)(e_{r+1}\dots
e_1)\dots(e_{r+s}\dots e_s)(e_{r+s+1}\dots e_{2+2s}).$$
If $r=s$ instead, the  computation goes through but the last
application of the lemma requires the second formula for the reduced
word. Hence the last factor needs to be replaced by the identity:
$$p_{r,2}p_{r+2,r}=\lambda^{2r+r^2}(e_r\dots e_1)(e_{r+1}\dots
e_1)(e_{r+2}\dots e_2)\dots(e_{2r}\dots e_r).$$
If $r>s$ we  apply the lemma iteratively, but
now
only  $s$
times (in spite of the $s+2$ factors) to
the first two
factors, the second and third factor and so on. We get
$$p_{r,2}p_{r+2,s}=\lambda^{2r+rs-2s}(e_r\dots e_1)\dots(e_{r+s-1}\dots 
e_s)
(e_{r+s}\dots e_{2s+1})(e_{r+s+1}\dots
e_{2+2s})=$$  
$$\lambda^{2(r-s)}p_{r,s}
(e_{r+s}\dots e_{1+2s})(e_{r+s+1}\dots
e_{2+2s})=p_{r,s}p^{(2s)}_{r-s,2},$$
and the formula is proved. Now  assume $r=s$. Then the right hand side of
the
desired formula reduces to $p_{r,r}$.
The same computation goes through except
at the last $s$-th iteration, where
$$p_{r,2}p_{r+2,r}=\lambda^{r^2}(e_r\dots e_1)\dots(e_{2r-1}\dots e_r)=
p_{r,r}.$$ 
The proof is now complete.
\medskip
It is now not difficult to interpret left tensoring by
the translates of $R$ in terms of the Jones tower.
\medskip

\noindent{\bf 5.8 Proposition} {\sl For $r,s\geq0$ and $\zeta\in 
M^{\otimes
r+s},$ we have:
\begin{description}
\item{\rm a)}
 for $r>s$,
$\Gamma_{r+2+s}(1_{M^{\otimes r}}\otimes R\otimes
1_{M^{\otimes s}}\circ\zeta)=
\lambda\Gamma_{r+s}(\zeta)p^{(2s)}_{r-s,2},$
\item{\rm b)} for $r=s$,
$\Gamma_{2r+2}(1_{M^{\otimes r}}\otimes R\otimes
1_{M^{\otimes r}}\circ\zeta)=\lambda\Gamma_{2r}(\zeta),$
\item{\rm c)} for $r<s$,
$\Gamma_{r+2+s}(1_{M^{\otimes r}}\otimes R\otimes
1_{M^{\otimes s}}\circ\zeta)=  
\lambda p^{(2r)}_{2,s-r}\Gamma_{r+s}(\zeta).$
\end{description}}\medskip
\noindent{\it Proof} For $r=0$ or $s=0$ the formula can be deduced
from $\Gamma_2(R)=\lambda$ and the computation of the tensor product 
given in Prop.\ 5.5. We can then assume $r,s>0$. In this case
a) and
b) follow from
the previous lemma,
 property $(5.12)$ and the fact that the $e_k$ commute with
$M_{s-1}$ for $k\geq 2s+1$ (recall $(3.1)$)
We prove c). Let us choose $\zeta$ of the form
$\zeta=\xi\otimes\eta_1\otimes\eta_2$ with $\xi,\eta_1\in M^{\otimes
r}$, $\eta_2\in M^{\otimes k}$, where $s=r+k$.
Then
$$\Gamma_{r+2+s}(1_{M^{\otimes r}}\otimes R\otimes
1_{M^{\otimes s}}\circ\zeta)=
\Gamma_{r+s+2}(\xi\otimes
R\otimes\eta_1\otimes\eta_2)=$$
$$\Gamma_{2r+2}(\xi\otimes
R\otimes\eta_1)p_{2r+2,k}\Gamma_k(\eta_2)=\lambda\Gamma_{2r}(\xi\otimes\eta_1)
p_{2r+2,k}\Gamma_k(\eta_2).$$
On the other hand
$$p_{2r+2,k}=\lambda^{k(2r+2)}((e_{2r+2}e_{2r+1})e_{2r}\dots
e_1)\dots((e_{2r+k+1}e_{2r+k})e_{2r+k-1}\dots e_{k})=$$
$$\lambda^{2k}
(e_{2r+2}e_{2r+1})(e_{2r+3}e_{2r+2})\dots(e_{r+s+1}e_{r+s})
p_{2r,k},$$
which implies
$$\Gamma_{r+2+s}(1_{M^{\otimes r}}\otimes R\otimes
1_{M^{\otimes s}}\circ\zeta)=$$
$$\lambda^{1+2(s-r)}\Gamma_{2r}(\xi\otimes\eta_1)
(e_{2r+2}e_{2r+1})(e_{2r+3}e_{2r+2})
\dots(e_{r+s+1}e_{r+s})
p_{2r,k}\Gamma_k(\eta_2)=$$
$$\lambda^{1+2(s-r)}(e_{2r+2}e_{2r+1})(e_{2r+3}e_{2r+2})
\dots(e_{r+s+1}e_{r+s})\Gamma_{r+s}(\zeta)=\lambda p^{(2r)}_{2,s-r}.$$
We next compute the operators in the Jones tower   corresponding to 
tensoring on the left by $1_{M^{\otimes r}}\otimes
R^*\otimes 1_{M^{\otimes s}}$.\medskip

\noindent{\bf 5.9 Proposition} {\sl For  $r,s\geq0$, $\zeta\in M^{\otimes
r+s+2}$ we have:
\begin{description}
\item{\rm a)} for $r>s$,
$\Gamma_{r+s}(1_{M^{\otimes r}}\otimes R^*\otimes
1_{M^{\otimes s}}\circ\zeta)=\lambda
E_{r+s}E_{r+s+1}(\Gamma_{r+s+2}(\zeta)(p^{(2s)}_{r-s,2})^*)$,
\item{\rm b)} for $r=s$, $\Gamma_{2r}(1_{M^{\otimes r}}\otimes R^*\otimes
1_{M^{\otimes r}}\circ\zeta)=\lambda 
E_{r+s}E_{r+s+1}(\Gamma_{2r+2}(\zeta))$,
\item{\rm c)} for $r<s$,
$\Gamma_{r+s}(1_{M^{\otimes r}}\otimes R^*\otimes
1_{M^{\otimes s}}\circ\zeta)=\lambda E_{r+s}E_{r+s+1}((p^{(2r)}_{2,s-r})^*
\Gamma_{r+s+2}(\zeta))$.
\end{description}}\medskip

\noindent{\it Proof } With respect to the inner products
the $N$--bimodule operators
$1_{M^{\otimes r}}\otimes R\otimes 1_{M^{\otimes s}}:
M^{\otimes r+s}\to M^{\otimes
r+s+2}$ and 
$1_{M^{\otimes r}}\otimes R^*\otimes 1_{M^{\otimes s}}:
M^{\otimes r+s+2}\to M^{\otimes
r+s}$ are adjoints of each other   Recall that,  by Prop. 10.5,
$\Gamma_p: M^{\otimes p}\to M_{p-1}$ is a $N$--bimodule unitary if
$M_{p-1}$ is regarded as a $N$--bimodule with inner product defined   
by the conditional expectation $E_{(p)}=EE_1\dots E_{p-1}: M_{p-1}\to N$.
Therefore we need to compute
$$\Gamma_{r+s} 1_{M^{\otimes r}}\otimes R^*\otimes
1_{M^{\otimes s}}\Gamma_{r+s+2}^*=
(\Gamma_{r+s+2}1_{M^{\otimes r}}\otimes
R\Gamma_{r+s}^*)^*.$$
Now by the previous lemma, for $r>s$,
$(\Gamma_{r+s+2}1_{M^{\otimes r}
}\otimes
R\Gamma_{r+s}^*): M_{r+s-1}\to M_{r+s+1}$ is right multiplication by
the element $A=\lambda p^{2s}_{r-s,2}\in M_{r+s+1}$.
Hence its adjoint, $r_A^*$,  is $$<X, r_A^*Y>=<XA,Y>=E_{(r+s+2)}(A^*X^*Y)=E_{(r+s+2)}(X^*YA^*)$$ 
where $X\in M_{r+s-1}$, $Y\in M_{r+s+1}$,
 as $A$ commutes with $N$. Hence
$$<X, r_A^*Y>=E_{(r+s)}(X^*E_{r+s}E_{r+s+1}(YA^*)),$$
as $X\in M_{r+s-1}$,
and this shows that
$ r_A^*Y= E_{r+s}E_{r+s+1}(YA^*)$.
The remaining cases follow similarly.\medskip

\noindent{\it Acknowledgements} We are grateful to L. Vainerman and J.--M. Vallin
for discussions. We also thank the referee for  asking whether   the noncommutative quotient  space
corresponding to index  $2$ subfactors could be identified with $D_1\backslash S_{-1}U(2)$ of \cite{Tomatsu}.
\medskip

\end{section}

\end{document}